\documentclass[a4paper,12pt]{article}

\usepackage{amsmath,amsfonts,graphicx}

\usepackage{float}

\numberwithin{equation}{section}
\numberwithin{table}{section}

\title{A Simple Method for High-Rank Families of Elliptic Curves with Specified Torsion}
\author{Allan J. MacLeod\\Statistics, O.R. and Mathematics Group,\\
University of the West of Scotland,\\High St.,  Paisley,\\Scotland.  PA1 2BE\\
(e-mail: allan.macleod@uws.ac.uk)}
\date{}

\begin{document}

\maketitle

\begin{abstract}
We describe a simple, but effective, method for deriving families of elliptic curves, with high rank, all of whose members have the same
torsion subgroup structure.
\end{abstract}

\newpage

\section{Introduction}

The standard elliptic curve over $\mathbb{Q}$ can be written
\begin{equation}
y^2+a_1\,x\,y+a_3y=x^3+a_2x^2+a_4x+a_6
\end{equation}
with $a_i \in \mathbb{Q}$. In fact, by suitable scaling, we could assume $a_i \in \mathbb{Z}$.

The rational points on this curve form an Abelian group, often denoted by $\Gamma$,
which is isomorphic to $T_G \oplus \mathbb{Z}^r$, where $r$ is the rank of the curve and $T_G$ is called the {\bf torsion subgroup} and is one of
\begin{itemize}
\item[{(i)}] $\hspace{1cm} \mathbb{Z}/n\mathbb{Z}, \hspace{0.5cm} n=1,2,\ldots,10,12$,
\item[{(ii)}] $\hspace{1cm} \mathbb{Z}/2\mathbb{Z} \oplus \mathbb{Z}/2n\mathbb{Z}, \hspace{0.5cm} n=1,2,3,4$.
\end{itemize}
as proven by Mazur \cite{mazur}, though seemingly first conjectured by Beppo Levi much earlier in the 20th century.

It is a standard conjecture that there exist elliptic curves with all possible
ranks. For specified torsion structures, however, the situation is less obvious.
The current situation for all $15$ torsion possibilities is summarized in Andrej Dujella's wonderful web-site at
\begin{center}
{\bf web.math.pmf.unizg.hr/~duje/tors/tors.html}
\end{center}
which should be referred to by anyone interested in high-rank elliptic curves.

As well as looking for specific curves of high rank, there is also the problem of trying to find infinite families of curves with high rank.
Many of the record families found have involved high-level mathematics, far above the level of a student or amateur.
It is the purpose of the current work to present an extremely simple method which can be understood by a much wider range of interested parties.

\section{Simple Elliptic Curve Forms}
Most discussions of the possible torsion structures start with Kubert's \cite{kub} list of parameterizations, which is based on the Tate normal form of an elliptic curve. The use of the Tate curve to
determine torsion parameterisations is partially described by Husem\"oller \cite{huse}.

The Tate form is
\begin{equation}\label{tate}
Y^2+(1-c)XY-bY=X^3-bX^2
\end{equation}
where $b,c$ are constant. This has the point $(0,0)$ which is assumed not of order $2$ or $3$.

If we define
\begin{equation*}
X=\frac{U}{4} \hspace{2cm} Y=\frac{4b+(c-1)U+V}{8}
\end{equation*}
the elliptic curve is in the symmetric form
\begin{equation}
V^2=U^3+((c-1)^2-4b)U^2+8b(c-1)U+16b^2
\end{equation}
with the point $(0,0)$ transformed to $(0,-4b)$.

Kubert provides formulae for $b$ and $c$ which give the torsion groups apart from $\mathbb{Z}/1\mathbb{Z}$, $\mathbb{Z}/2\mathbb{Z}$ and
$\mathbb{Z}/3\mathbb{Z}$. Unfortunately, these can lead to formulae which are more complex than need be.

In this section, we discuss providing simple forms which assist in developing families with rank greater than $0$. We use $\mathbb{Z}/8\mathbb{Z}$
as an example. The point at infinity is the identity, and there is $1$ point of order two, $2$ points of order four and $4$ points
of order eight.

The presence of a single point of order $2$ means that we can assume, without loss of generality, that the curve has form
\begin{equation}
y^2=x^3+Ax^2+Bx
\end{equation}
with $A,B \in \mathbb{Z}$. As such, the Nagell-Lutz theorem \cite{siltate} tells us that all the torsion points have integer coordinates.
We also have that, if $P=(Z,W)$ is a rational point on the curve, then the x-coordinate of $2P$ is
\begin{equation}
\frac{(Z^2-B)^2}{4W^2}
\end{equation}
and thus is a rational square.

If $P$ is a point of order $4$, then $x(2P)=0$ so $Z^2=B$, and thus $B$ must be an integer square. Thus we can redefine the curve as
\begin{equation}
y^2=x^3+Ax^2+B^2x
\end{equation}
and, substituting $x=B$ gives $Q^2=B^2(A+2B)$. Thus, there must exist $C \in \mathbb{Z}$ with $A=C^2-2B$, so the curve is now
\begin{equation}
y^2=x^3+(C^2-2B)x^2+B^2x
\end{equation}

Since twice a point of order $8$ gives a point of order $4$, using the same logic about doubling shows that $B=D^2$ with $D \in \mathbb{Z}$,
and the curve is thus
\begin{equation}
y^2=x^3+(C^2-2D^2)x^2+D^4x
\end{equation}

If $(F,G)$ is a point of order $8$ then $F$ satisfies the quartic
\begin{equation*}
x^4-4D^2x^3+2D^2(3D^2-2C^2)x^2-4D^6x+D^8=0
\end{equation*}
which factors into
\begin{equation}
(x^2-2D(C+D)x+D^4)(x^2+2D(C-D)x+D^4)=0
\end{equation}

The first quadratic gives a rational solution iff the discriminant is a rational square, giving
\begin{equation*}
C^2+2CD=\Box
\end{equation*}

Consider the related quadric $t^2+2t=s^2$, which contains $(0,0)$. The line through the origin $s=kt$ meets the curve again at
\begin{equation}
t=\frac{2}{1-k^2}=\frac{2m^2}{m^2-n^2} \hspace{2cm} k=\frac{n}{m}
\end{equation}
so, taking $C=2m^2$ and $D=m^2-n^2$ gives
\begin{equation}
y^2=x^3+2(m^4+2m^2n^2-n^4)x^2+(m^2-n^2)^4x
\end{equation}

Finally, set $r=m/n$, $x=x/n^4$ and $y=y/n^6$ to give
\begin{equation}\label{z8ec}
y^2=x^3+2(r^4+2r^2-1)x^2+(r^2-1)^4x \hspace{2cm} r \neq 0, \pm 1
\end{equation}
which has a point of order $2$ at $x=0$, two points of order $4$ at $x=(r^2-1)^2$ and four points of order $8$
at $x=-(r-1)(r+1)^3$ and $x=-(r+1)(r-1)^3$.

\section{Basic Method}

In May $2013$, Dujella and Juan Carlos Peral submitted a very interesting preprint \cite{dujper} to the arXiv server. It concerned finding families of curves with
$\mathbb{Z}/8\mathbb{Z}$ or $\mathbb{Z}/2\mathbb{Z} \oplus \mathbb{Z}/6\mathbb{Z}$ torsion and rank two.

Unfortunately, they gave many formulae but little information on the precise nature of their methodology. The present work was inspired
by the statement "By searching on several homogeneous spaces ...".
This links the method to the standard method of 2-descent for an individual curve, as described by Cremona \cite{crem}.

For an equation of the form
\begin{equation*}
y^2=x^3+Ax^2+Bx \hspace{2cm} A,B \in \mathbb{Z}
\end{equation*}
it is known that the $(x,y)$ coordinates of a rational point are of the form
\begin{equation}
x=\frac{du^2}{v^2} \hspace{2cm} y=\frac{duw}{v^3}
\end{equation}
where $d,u,v,w \in \mathbb{Z}$ and $d$ is squarefree, and
\begin{equation*}
\gcd(d,u)=\gcd(d,v)=\gcd(u,v)=1
\end{equation*}
Some simple algebra shows that this implies $d | B$.

This suggests searching for points on \eqref{z8ec} of the form
\begin{equation}
x=\frac{d(r)\,u^2(r)}{v^2(r)}
\end{equation}
where $u(r)$ and $v(r)$ are simple polynomials and $d(r)$ is of the form
\begin{equation}
d(r)= \pm 2^{e} (r-1)^{f} (r+1)^g
\end{equation}
with $e,f,g \in \mathbb{Z}$.

Since many curves with this torsion have rank $0$, we are very unlikely to find such a point directly. We can, however, look for such
points which, when substituted into the right-hand-side of \eqref{z8ec}, it simplifies to
\begin{equation}
F^2(r) (a r^2+br+c)
\end{equation}
with $F(r)$ being a rational function and $a,b,c \in \mathbb{Q}$.

Thus we will have a point of infinite order (apart from a small number of $r$ values) when
\begin{equation}
ar^2+br+c=\Box=t^2
\end{equation}
and we attempt to parameterize such quadratics. Note that we do not always have a quadratic which is soluble.

Suppose we find a solution $r=p, t=q$, then the line $t=q+k(r-p)$ meets the curve at one further point giving
\begin{equation*}
r=\frac{p\,k^2-2\,q\,k+a\,p+b}{k^2-a}
\end{equation*}

As an example, we find quickly that $x=1-r^2$ gives
\begin{equation}
y^2=r^4(r+1)^2(r-1)^2(5-r^2)
\end{equation}
so a rational point occurs when $5-r^2=d^2$.

There is a solution $r=1,\, d=2$, and so the line $d=2+t(r-1)$ meets the curve again where
\begin{equation}
r=\frac{t^2-4t-1}{t^2+1}
\end{equation}

Substituting into \eqref{z8ec} and defining
\begin{equation*}
x=\frac{u}{(t^2+1)^4} \hspace{2cm} y=\frac{v}{(t^2+1)^6}
\end{equation*}
gives the elliptic curve
\begin{equation*}
v^2=u^3+4(t^8-16t^7+60t^6-112t^5+62t^4+112t^3+60t^2+16t+1)u^2+
\end{equation*}
\begin{equation}\label{ecz8r1}
256t^4(t-2)^4(2t+1)^4u
\end{equation}

\begin{table}
\begin{center}
\caption{$Z8$ points and quadratics}
\begin{tabular}{lr}
$\,$&$\,$\\
Point&Quadratic\\
$\,$&$\,$\\
$-(r^2-1)^2$&$2r^2-1$\\
$\,$&$\,$\\
$-(r^2-1)^3$&$5-r^2$\\
$\,$&$\,$\\
$(r+1)^4$&$r^2+4$\\
$\,$&$\,$\\
$\frac{-(4r^2-3)^2}{16}$&$32r^2-25$\\
$\,$&$\,$\\
$\frac{-(r^2-1)(r-1)^2(r+2)^2}{(r-2)^2}$&$16-7r^2$\\
$\,$&$\,$\\
$(r+1)^3(1-3r)$&$-3r^2-14r+5$\\
$\,$&$\,$\\
$-(r+1)^2(r^2+2r-1)$&$-2r^2-4r+2$\\
$\,$&$\,$\\
$1-r^4$&$2r^2+2$\\
$\,$&$\,$\\
$-(r+1)^2(r^2+r-1)$&$-r^2-r+1$\\
$\,$&$\,$\\
$\frac{27(1-r^2)}{2}$&$150-6r^2$\\
$\,$&$\,$\\
$\frac{(r-1)^4(2r+1)}{1-2r}$&$1-4r^2$\\
$\,$&$\,$\\
$\frac{2(r^2-1)^4}{2-r^2}$&$4-2r^2$
\end{tabular}
\end{center}
\end{table}

From the final term, it is clear that we cannot have $t=0,2,-1/2$ as they give a singular curve. The point of order $4$ becomes
$16t^2(t-2)^2(2t+1)^2$ and the two points of order $8$ are now $16t^3(t-2)^3(2t+1)$ and $16t(t-2)(2t+1)^3$. The point $1-r^2$ is
$4t(t-2)(t^2+1)^2(2t+1)$. We have to ensure this is not equal to one of the torsion points and this condition also excludes
$t=\pm 1, -3, 1/3$.

Table $3.1$ gives the points found and the soluble quadratics. From \eqref{z8ec}, it is clear that we can ignore the variations in signs
for $r$.

For each of these, it is straightforward to find a parametrization of the quadratic and hence a form of the elliptic curve
similar to \eqref{ecz8r1}.

\section{Rank $2$ Curves}
Consideration of equation \eqref{ecz8r1} shows that we can repeat the method of the previous section on this and similar curves. We will have
more complicated functions to factorize but the essential idea goes through.

For the rank $1$ curves found, we quickly found that the fifth point in the Table gave a quadratic term. We summarize the relevant details.
The quadratic $16-7r^2=\square$ can be parametrized by
\begin{equation*}
r=\frac{k^2-6k-7}{k^2+7}
\end{equation*}
using the solution $r=1, \, d=3$.

Substituting into \eqref{z8ec} and clearing denominators gives
\begin{equation*}
v^2=u^3+4(k^8-24k^7+116k^6-264k^5-458k^4+1848k^3+
\end{equation*}
\begin{equation}\label{ecz8r15}
5684k^2+8232k+2401)u^2+256k^4(k-3)^4(3k+7)^4u
\end{equation}

We applied the method to this curve with
\begin{equation}
d = \pm 2^e \, k^f \, (k-3)^g \, (3k+7)^h
\end{equation}
for small integer $e,f,g,h$.

For $u=k^2(k+1)^2(k-3)(3k+7)^2$, we have a rational point on \eqref{ecz8r15} if $4k^2+9k+1=d^2$. From the solution $k=0, \, d=1$, we have
\begin{equation*}
k=\frac{9-2s}{s^2-4}
\end{equation*}
which gives the following form of the curve
\begin{equation}\label{ecz8r2}
Z^2=W^3+4(2401s^{16}-16464s^{15}+19992s^{14}+241584s^{13}-891820s^{12}-
\end{equation}
\begin{equation*}
1083024s^{11}+7095080s^{10}+5339280s^9-29327658s^8-63355440s^7+
\end{equation*}
\begin{equation*}
162576392s^6+280936080s^5-869091628s^4+139426320s^3+
\end{equation*}
\begin{equation*}
1538092920s^2-2179390416s+1166778337)W^2+
\end{equation*}
\begin{equation*}
256(s-1)^4(s+2)^4(s-2)^4(s+3)^4(2s-9)^4(3s-7)^4(7s+1)^4W
\end{equation*}
together with the two points of infinite order
\begin{equation}
W_1=\frac{16\,C\,(7s^4+12s^3-98s^2-156s+571)^2}{(21s^4-12s^3-110s^2+12s+201)^2}
\end{equation}
where $C=(s+3)(s-1)^3(s+2)^3(s-2)^3(2s-9)(3s-7)(7s+1)^3$
and
\begin{equation}
W_2=(s+2)(2-s)(s+3)(s-1)^2(s^2-2s+5)^2(2s-9)^2(3s-7)(7s+1)^2
\end{equation}

Investigations of the torsion points and other considerations show that $s$ cannot be in the set $\{1, \pm2, -3, 9/2, 7/3, -1/7 \}$.

Specializing with $s=4$, gives the curve
\begin{equation*}
Z^2=W^3+5379102933124W^2+456366899570319360000W
\end{equation*}
and
\begin{equation*}
W_1=\frac{-561850417010211840}{9591409} \hspace{2cm} W_2=-537247620
\end{equation*}

Transforming to the minimal model of this elliptic curve, we find the height matrix has determinant $90.592$ showing the rank is at least $2$.

After this rank $2$ family was found, the author informed Prof. Dujella, who sent the e-mail to Dr. Peral. He quickly discovered only the second
known example of a curve with $\mathbb{Z}/8\mathbb{Z}$ torsion and rank $6$. This came some $7$ years after Noam Elkies discovered the first example!

The above rank $2$ family was found early in this investigation. Over time, the methodology and software were developed and, more recently, a second rank $2$
family was discovered. We will just give the relevant details - the interested reader just has to apply the same ideas as above.

We have
\begin{equation*}
x=\frac{(r+1)^2(2r-1)^2(5r-4)}{4(3r-1)}
\end{equation*}

\begin{equation*}
r=\frac{17k^2-60k}{15k^2-225}
\end{equation*}

\begin{equation*}
u=-9(2k-5)^2(4k-15)^3(8k+15)^3
\end{equation*}

\begin{equation*}
k=\frac{-5s^2-80s+580}{2s^2+184}
\end{equation*}

\section{Curves with $\mathbb{Z}/2\mathbb{Z} \oplus \mathbb{Z}/6\mathbb{Z}$ torsion}
For curves with $\mathbb{Z}/2\mathbb{Z} \oplus \mathbb{Z}/6\mathbb{Z}$ torsion, the methodology is the same as before, so only the relevant results are quoted.

Kubert's formulae for this torsion are
\begin{equation*}
b=c+c^2 \hspace{2cm} c=\frac{10-2\alpha}{\alpha^2-9}
\end{equation*}
which, after some manipulation, gives the curve
\begin{equation*}
g^2=h^3-2(\alpha^4-12\alpha^3+30\alpha^2+36\alpha-183)h^2+(\alpha+3)(\alpha-9)(\alpha-1)^3(\alpha-5)^3h
\end{equation*}

Numerical experimentation suggested a simpler form could be found, which is
\begin{equation}\label{z26ec}
y^2=x^3-2(r^4-6r^2-3)x^2+(r^2-1)^3(r^2-9)x
\end{equation}
where $\alpha=2r+3$.

With this form, we have found over $50$ different x-points which leave a quadratic to be a square. The software
then produced several rank-2 solutions, many of which were essentially equivalent. There are $7$ solutions which seem to be independent,
and are given in a similar form to the description in the previous sections:

Solution 1:

\begin{equation*}
x=\frac{(r-3)(r-1)^2(r+1)(r+2)}{r-2} \hspace{2cm}  r=\frac{k^2-6k-3}{k^2+3},
\end{equation*}

\begin{equation*}
u=-4k^2(k^2-9)^2(2k^2-3k+3) \hspace{2cm}  k=\frac{3s^2-12s+9}{s^2-2}
\end{equation*}

\vspace{1cm}

Solution 2:

\begin{equation*}
x=(r^2-1)^3/r^2 \hspace{2cm}  r=\frac{k^2+9}{27-3k^2}  ,
\end{equation*}

\begin{equation*}
u=144(k^2-18)(4k^2-45)(5k^2-36) \hspace{2cm}  k=\frac{3s^2-6s+24}{s^2-8}
\end{equation*}

\vspace{1cm}

Solution 3:

\begin{equation*}
x= \frac{(r-7)(r-3)(r^2-1)^2}{(r-5)^2} \hspace{2cm}  r= \frac{k^2-12k-48}{k^2-6},
\end{equation*}

\begin{equation*}
u= -16(k-9)(k+3)(2k+7)(k^2+6k+15)(2k^2-6k-33)  \hspace{1cm}  k=  \frac{3s^2+2s+7}{1-s^2}
\end{equation*}

\vspace{1cm}

Solution 4:

\begin{equation*}
x= (r+1)^2(r+3)(r-0.5),   \hspace{2cm}  r=\frac{k^2+54}{108-3k^2}    ,
\end{equation*}

\begin{equation*}
u= 3600(k^2-81)(k^2-27)(4k^2-189)  \hspace{2cm}  k= \frac{6s^2-30s+240}{s^2-40}
\end{equation*}

\vspace{1cm}

Solution 5:

\begin{equation*}
x=  \frac{2(r^2-9)(r^2-1)(r+1)}{3r-6}       \hspace{2cm}  r= \frac{k^2-12k-48}{k^2-6},
\end{equation*}

\begin{equation*}
u= -4(k-9)^2(k+3)^2(k+6)^2(2k^2-6k-33)  \hspace{1cm}  k= \frac{3s^2+6s}{2-s^2}
\end{equation*}

\vspace{1cm}

Solution 6:

\begin{equation*}
x=\frac{(r-3)(r-1)(r^2-1)^2}{(r-2)^2}      \hspace{2cm}  r= \frac{5k^2-15}{k^2-9}
\end{equation*}

\begin{equation*}
u= 2160(k^2-4)^2(4k^2-21)  \hspace{2cm}  k= \frac{2s^2-10s+80}{s^2-40}
\end{equation*}

\vspace{1cm}

Solution 7:

\begin{equation*}
x=\frac{2(r+1)(r-3)(r-1)^3}{x-2}     \hspace{2cm}  r=\frac{5k^2-4}{k^2-2}
\end{equation*}

\begin{equation*}
u= 16(k^4-1)(2k^2-1)^3 \hspace{2cm}  k= \frac{s^2+2s+8}{8-s^2}
\end{equation*}

\vspace{1cm}

\section{Curves with $\mathbb{Z}/7\mathbb{Z}$ torsion}
Analysis of Dujella's tables of ranks shows that this torsion forms a trio with the previous two torsion subgroups in that they have a similar structure and results. We have been
unable to find a simpler form than that given by Kubert from the Tate form \eqref{tate}.

Set $b=r^3-r^2$ and $c=r^2-r$, then we can transform to the symmetric form
\begin{equation}\label{ecz7}
y^2=x^3+(r^4-6r^3+3r^2+2r+1)x^2+8r^2(r-1)(r^2-r-1)x+16r^4(r-1)^2
\end{equation}

For curves with no 2-torsion, rational points do not have such a structured form as before. We thus choose $d(r)$ from divisors of
\begin{equation*}
16r^4(r-1)^2  \hspace{1cm} \mbox{or} \hspace{1cm} 16r^4(r-1)^2(r^2-r-1)
\end{equation*}

Table $6.1$ list the rational points $x$ and the related quadratic to be made square.
\begin{table}
\begin{center}
\caption{$\mathbb{Z}/7\mathbb{Z}$ points giving rank $1$ families}
\begin{tabular}{lr}
$\,$&$\,$\\
Point&\hspace{1cm}Quadratic\\
$\,$&$\,$\\
$4r - 4$&$4r^2-3$\\
$-r^2$&$r^2-12r+9$\\
$-4r^2$&$r^2-8r+4$\\
$3r^2 - 2r - 1$&$9r^2-6r-2$\\
$-r^4 + 6r^3 - 5r^2$&$-2r^2+10r+1$\\
$-4r^4 + 12r^3 - 8r^2$&$-3r^2+6r+1$\\
$4r^2(r-1)/(2r - 1)^2$&$100r^2-116r+225$\\
$16r^2(r-1)/(r + 1)^2$&$25r^2+66r+9$\\
$4r^2(r-1)/(r^2 - r + 1)$&$r^2-r+1$\\
$4r^2(r-1)/(r^2-9r+9)$&$r^2-9r+9$\\
$-4r^2(r-1)(3r-4)/(r-2)^2$&$9r^2-84r+100$
\end{tabular}
\end{center}
\end{table}

Unfortunately, despite extensive computation, we have been unable to find a rank $2$ family with this torsion. We have, however, found points
which give curves of rank $2$ subject to a quartic being made square. These quartics are all equivalent to elliptic curves with rank usually $1$, but
sometimes $2$.

In \cite{lec}, Odile Lecacheux found parametric forms for $r$ and points $x_1$ giving rank-one families, using fibrations of related surfaces. Possibly the simplest of these gave
\begin{equation*}
r=\frac{k^2-1}{k(k-2)} \hspace{2cm} x_1=\frac{-4(k-1)(2k-1)}{k^2(k-2)}
\end{equation*}
Again these points did not yield a family with rank $2$.

\section{Other Torsion Subgroups}
The curves with torsion subgroup isomorphic to $Z/5Z$, $Z/6Z$, or $Z/2Z \times Z/4Z$ also form a related group. The ranks are higher
so we have to perform the basic strategy more often, but we have found rank $3$ families for all three torsion subgroups. Again, we
just give the relevant details.

\subsection{$Z/5Z$}
The parametric form is just that of Kubert and is
\begin{equation}
y^2=x^3+(r^2-6r+1)x^2+8r(r-1)x+16r^2
\end{equation}

There are 2 independent solutions that have been found.

Solution $1$:

\begin{equation*}
x=\frac{-12r}{r-3} \hspace{2cm} r=\frac{3*k^2-12}{k^2-1},
\end{equation*}

\begin{equation*}
u=4(k^2-4)(k-1)^2 \hspace{2cm} k=\frac{(s-1)^2}{s^2-17},
\end{equation*}

\begin{equation*}
h=-16(s-5)(s+3)(s+7)(3s-11)(s^2-s-8)^2 \hspace{1cm} s=\frac{5t^2-40t+131}{t^2-33}
\end{equation*}
where the third elliptic curve is assumed to be of the form $g^2=$ cubic in $h$, and that denominators have been cleared.

\vspace{1cm}

Solution $2$:

\begin{equation*}
x=\frac{4r}{1-3r)} \hspace{2cm}r=\frac{k^2-144}{3(k^2-36)},
\end{equation*}

\begin{equation*}
u=4(k-6)^2(k^2-144)   \hspace{2cm} k=\frac{12(s^2-18s+81)}{s^2-41},
\end{equation*}

\begin{equation*}
h=-82944(s-4)(s-5)(9s-61)(3s^2-36s+121)^2 \hspace{1cm} s=\frac{t^2-412t-2213}{3t^2-3}
\end{equation*}

\subsection{$Z/6Z$}
For this torsion, we transform the form derived from Kubert to the much simpler
\begin{equation}
y^2=x^3+(r^2-3)x^2+(3-2r)x
\end{equation}

We have found only one independent family of rank $3$, though it appears in many forms. The relevant formulae are

\begin{equation*}
x= 4-3r \hspace{2cm}r= \frac{k^2-4k-24}{k^2-9} ,
\end{equation*}

\begin{equation*}
u= 9(k-3)(k+3)^3     \hspace{2cm} k= \frac{2s^2+6s+12}{4-s^2}
\end{equation*}

\begin{equation*}
h= (s^2-6s-24)^2(9s^4-24s^3-116s^2+336s+864) \hspace{1cm} s=  \frac{-t^2-58t+347}{t^2+47}
\end{equation*}

\subsection{$Z/2Z \times Z/4Z$}
The parametric form used is
\begin{equation}
y^2=x^3+(r^2+1)x^2+r^2x=x(x+1)(x+r^2)
\end{equation}

We found $4$ solutions, but three are essentially the same, so there are only $2$ independent solutions.

Solution $1$:

\begin{equation*}
x=2r+1, \hspace{2cm} r= \frac{2-k^2}{k^2-4} ,
\end{equation*}

\begin{equation*}
u=-(k+2)^4 \hspace{2cm} k= \frac{-6}{s^2+4},
\end{equation*}

\begin{equation*}
h= 8(s^2+1)^2(s^4+8s^2-2) \hspace{1cm} s= \frac{ t^2-6t+1}{t^2-1}
\end{equation*}

\vspace{1cm}

Solution $2$:

\begin{equation*}
x= -3r  \hspace{2cm}r=\frac{6k+30}{9-k^2} ,
\end{equation*}

\begin{equation*}
u= 3(k-11)(k+3)^2(k+5)  \hspace{2cm} k= \frac{s^2+11}{1-s^2},
\end{equation*}

\begin{equation*}
h= -64(s^2-4)^2(s^2-1)(s^2+2)^2  \hspace{1cm} s= \frac{t^2-4t-1}{t^2+1}
\end{equation*}


\newpage

\begin{thebibliography}{199}

\bibitem{crem} Cremona J.E., {\it Algorithms for Modular Elliptic Curves}, Cambridge University Press, $1997$.

\bibitem{dujper} Dujella A. and Peral J.C., {\it Elliptic curves with torsion group $Z/8Z$ or $Z/2Z \times Z/6Z$}, arXiv preprint $1306:0027$ (2013).

\bibitem{huse} Husem\"oller D., {\it Elliptic Curves}, Springer-Verlag, $2004$.

\bibitem{kub} Kubert D.S., {\it Universal bounds on the torsion of elliptic curves}, Proc. London Math. Soc., {\bf 33} (1976) 193-237.

\bibitem{lec} Lecacheux O., {\it Rang de familles de courbes elliptiques}, Acta Arith. {\bf 109} (2003) 131-142.

\bibitem{mazur} Mazur B., {\it Modular curves and the Eisenstein ideal}, Pub. Math. IHES., {\bf 47} (1977) 33-186.

\bibitem{siltate} J.H. Silverman and J. Tate, {\it Rational Points on Elliptic Curves}, Springer, $1992$.


\end{thebibliography}
\end{document}